\documentclass[article,reqno,11pt]{amsart}

\usepackage{amssymb,amsfonts,amsmath,amsthm}
\usepackage[all,arc]{xy}
\usepackage{enumerate}
\usepackage{mathrsfs}
\usepackage[toc,page]{appendix}
\usepackage[left=3cm, right=3cm, bottom=3cm]{geometry}
\usepackage{tabularx}
\usepackage{url}
\usepackage{color}
\usepackage{tikz-cd}

\newtheorem{thm}{Theorem}[section]

\newtheorem*{thm*}{Theorem}
\newtheorem{cor}[thm]{Corollary}

\newtheorem{lem}[thm]{Lemma}

\theoremstyle{definition}
\newtheorem{defn}[thm]{Definition}

\newtheorem{con}[thm]{Construction}

\theoremstyle{remark}
\newtheorem{rem}[thm]{Remark}

\newtheorem*{idea*}{Idea}

\makeatletter
\let\c@equation\c@thm
\makeatother
\numberwithin{thm}{section}
\numberwithin{equation}{section}
\title{\textsc{Degree of W-operator and Noncrossing Partition}}

\author{Hao Sun}

\begin{document}
\maketitle
\flushbottom

\begin{abstract}
	Goulden and Jackson first introduced the cut-and-join operator. The cut-and-join is widely used in studying the Hurwitz number and many other topological recursion problems. Mironov, Morosov and Natanzon give a more general construction and call it $W$-operator $W([n])$. As a special case, the cut-and-join operator is $W([2])$. In this paper, we study the structure of $W([n])$. We prove that $W([n])$ can be written as the sum of $n!$ terms and each term corresponds uniquely to a permutation in $S_n$. We also define the degree of each term. We prove that there is a correspondence between the terms of $W([n])$ with highest degree and the noncrossing partitions.
\end{abstract}


\section{Introduction}

The cut-and-join operator $\Delta$
\begin{align*}
	\Delta=\frac{1}{2}\sum_{i,j \geq 1}\left((i+j)p_i p_j\frac{\partial}{\partial p_{i+j}}+ ijp_{i+j}\frac{\partial^2}{\partial p_i \partial p_j}\right)
\end{align*}
is introduced by Goulden and Jackson \cite{MR1249468,GJ4}. It is an infinite sum of differential operators in variables $p_i$, $i \geq 1$. This operator plays an important role in calculating the simple Hurwitz number \cite{GJ5,GJ6,GJ7,MR1396978} and many other enumerative geometry problems \cite{GJV,Lando,Sun1}.

Mironov, Morosov and Natanzon \cite{MR2864467,Mir} construct the $W$-operators $W(\lambda)$, where $\lambda=(1^{i_1}2^{i_2}\dots n^{i_n})$ is a partition of a positive integer $n$. The $W$-operators are differential operators acting on the formal power series $\mathbb{C}[[X_{ij}]]_{i,j \geq 1}$, where $X_{ij}$ are coordinate functions on the positive-half-infinite matrix. A subring of $\mathbb{C}[[X_{ij}]]_{i,j \geq 1}$ is $\mathbb{C}[p_1,p_2,...]$, where $p_k=Tr(X^k)$ and $X=(X_{ij})_{i,j \geq 1}$. Let $[n]=(n^1)$ be the partition of $n$. A direct calculation shows that $W([2])$ is the cut-and-join operator $\Delta$ on the ring $\mathbb{C}[p_1,p_2,...]$. In this paper, we study the structure of the operators $W([n])$, $n \geq 1$. We consider $W([n])$ as operators on the ring $\mathbb{C}[p_1,p_2,...]$.

In \S 2, we review the natural quiver structure of permutations. Based on the quiver structure, we give a construction (Construction \ref{202}) that all permutations in $S_{n+1}$ can be constructed from permutations in $S_n$. Generally speaking, given any permutation in $S_n$, we construct $n+1$ permutations in $S_{n+1}$. The key part in \S 2 is Construction \ref{202}. Construction \ref{202} plays an important role in proving the structure theorem of $W([n])$ (Theorem \ref{307}).

In \S 3, we review the definition and some properties of $W([n])$. We prove the following theorem about the structure of $W([n])$.
\begin{thm*}{\textnormal{\textbf{\ref{307}}}}
	$W([n])$ is a well-defined operator on $\mathbb{C}[p_{1},p_{2},\dots]$. It can be written as the sum of $n!$ summations $FS_\alpha$, each of which corresponds to a unique quiver $\hat{Q}_\alpha$, or equivalently a unique permutation $\alpha \in S_n$.
\end{thm*}

In \S 4, we define the degree of each summation (term) of $W([n])$. We show that the maximal degree is $n+1$. The maximal degree summations are used to study the generating function of minimal $d$-Hurwitz number \cite{MR1797682}, where $d$ stands for $d$-cycles.

In \S 5, we study the number of summations of $W([n])$ with maximal degree. We show that there is an one-to-one correspondence between summations with maximal degree and noncrossing permutations of $[n]=\{1,\dots,n\}$ (Theorem \ref{504}). By Simion's study of noncrossing partition \cite{Sim}, we show that the number of summations with maximal degree is exactly the Catalan number (Corollary \ref{505}).

\section{Permutation Group and Quiver}
In this section, we first review the quiver structure of permutations. Then, we give a new construction about all permutations in $S_{n}$ by induction. Roughly speaking, we construct $n$ distinct permutations in $S_n$ from a given permutation in $S_{n-1}$. This construction comes from some operations on the quiver, and the idea of this construction comes from the calculation of $W$-operator (Remark \ref{306}). We will discuss the idea in the next section.

We begin with the definition of the quiver. A \emph{quiver} is a directed graph. So, as usual, a \emph{quiver} $Q=(V,A,s,t)$ is a quadruple, where $V$ is the set of vertices, $A$ is the set of arrows, $s$ and $t$ are two maps $A \to V$. If $a \in A$, $s(a)$ is the source of this arrow and $t(a)$ is the target. We assume $V$ and $A$ to be finite sets. If $B$ is a subset of $A$, $V_{B}=\{s(a),t(a),a \in B\}$, then we call $(V_{B},B,s',t')$ the \emph{subquiver} of $Q$, where $s'=s|_B$, $t'=t|_B$. A quiver $Q=(V,A,s,t)$ is \emph{connected} if the underlying undirected graph of $Q$ is connected. A connected quiver $Q=(V,A,s,t)$ is a \emph{loop}, if for any vertex $v \in V$, there is a unique arrow $a \in A$ such that $s(a)=v$ and a unique arrow $b \in A$ such that $t(b)=v$. A \emph{chain} is obtained by omitting a single arrow in a loop.

Any permutation $\alpha \in S_n$ has a natural quiver structure $Q_\alpha=(V_\alpha,A_\alpha,s,t)$, where $V_\alpha=\{1,\dots,n\}$ and $A_\alpha=\{i \rightarrow \alpha(i), 1 \leq i \leq n \}$. Notice that the quiver $Q_\alpha$ only contains loops. Since we want to use induction, we construct another quiver $\hat{Q}_\alpha$ from $Q_\alpha$.
\begin{con}\label{201}
Given $\alpha \in S_{n}$, let $Q_\alpha$ be the corresponding quiver. There is a unique arrow $a$ in $Q_{\alpha}$ such that $s(a)=1$. We substitute this arrow by a new one $\hat{a}$, where $s(\hat{a})=n+1$ and $t(\hat{a})=t(a)$. Denote by $\hat{A}_\alpha$ the new set of arrows. The new vertex set is $\hat{V}_\alpha=\{1,\dots,n,n+1\}$. Denote by $\hat{Q}_{\alpha}$ the new quiver, $\hat{Q}_{\alpha}=(\hat{V}_{\alpha}, \hat{A}_{\alpha},s,t)$.
\end{con}
We give an example about this construction. Take $\alpha =(3 2 1) \in S_{3}$, then $Q_{\alpha}$ is
\[
\begin{tikzcd}
Q_\alpha : \quad 3 \arrow[r] & 2 \arrow[r] & 1 \arrow[ll,bend right] \text{ }.
\end{tikzcd}
\]
$\hat{Q}_{\alpha}$ is
\[
	\begin{tikzcd}
	\hat{Q}_{\alpha} : \quad 4 \arrow[r] & 3 \arrow[r] & 2 \arrow[r] & 1 \text{ } .
	\end{tikzcd}
\]
Clearly, $Q_\alpha$ is a loop and $\hat{Q}_{\alpha}$ is a chain. In general, $\hat{Q}_\alpha$ consists of one chain and possibly a number of loops. The chain in $\hat{Q}_\alpha$ always start from $n+1$ and stop at $1$. Clearly, we can construct $Q_\alpha$ uniquely from $\hat{Q}_\alpha$. Thus we have an one-to-one correspondence between permutations $\alpha$ and quivers $\hat{Q}_{\alpha}$.

Now we will work on the quiver $\hat{Q}_\alpha$ and construct $n+1$ quivers corresponding to permutations in $S_{n+1}$. Generally speaking, if we want to construct a new quiver $\hat{Q}_\beta$ representing an element $\beta \in S_{n+1}$, we should add one more vertex $n+2$ into $\hat{V}_\alpha$ and add arrows $a_1,a_2$ in $\hat{A}_\alpha$ such that
\begin{align*}
s(a_1)=n+2, \quad t(a_2)=n+1,
\end{align*}
where $a_1,a_2$ can be the same arrow. Here is the construction.
\begin{con}\label{202}
Given any $\alpha \in S_{n}$, we write $\alpha$ as the product of disjoint cycles $\alpha=\alpha_{1}\alpha_{2}\dots\alpha_{k}$. We assume $1 \in \alpha_{1}$. So, the corresponding subquiver for $\alpha_{1}$ in $\hat{Q}_{\alpha}$ is the chain as following
\[
\begin{tikzcd}
\hat{Q}_{\alpha_1}: \quad n+1 \arrow[r] &  \cdots \arrow[r]  & 1  \text{ }.
\end{tikzcd}
\]
\begin{itemize}
\item{\textbf{Case 0}}

We extend the quiver for $\alpha_{1}$ directly
\[
\begin{tikzcd}
\hat{Q}_{\beta_1} : \quad n+2 \arrow[r] & n+1 \arrow[r] &  \cdots \arrow[r]  & 1  \text{ }.
\end{tikzcd}
\]

Clearly, this subquiver represents a well-defined cycle $\beta_1$ (by replacing $n+2$ by $1$). In this way, we construct a permutation $\beta \in S_{n+1}$, where $\beta=\beta_1\alpha_2\dots\alpha_k$. In this case, $a_1, a_2$ are the same arrow
\[
\begin{tikzcd}
a_1=a_2 : & n+2 \arrow[r] &  n+1 \text{ }.
\end{tikzcd}
\]
\end{itemize}
Next we consider the general case. The idea is cutting an arrow in $\hat{Q}_\alpha$ and reconnect the chain and loops in $\hat{Q}_{\alpha}$. Since there are $n$ choices of arrows in $\hat{Q}_\alpha$, we can construct $n$ permutations.

We choose an arbitrary arrow $a:i \to j$ in $\hat{Q}_{\alpha}$.
\begin{itemize}
\item{\textbf{Case 1} (\emph{Cut Case}), $a \in \hat{Q}_{\alpha_{1}}$}

In this case, $\hat{Q}_{\alpha_{1}}$ is
\[
\begin{tikzcd}
\hat{Q}_{\alpha_1} : \quad n+1 \arrow[r] &  \cdots \arrow[r]  & i \arrow[r] & j \arrow[r] & \cdots \arrow[r] & 1 \text{ }.
\end{tikzcd}
\]
 	
First, cut the arrow $i \to j$. We get
\[
\begin{tikzcd}
n+1 \arrow[r] &  \cdots \arrow[r]  & i \text{ }, & j \arrow[r] & \cdots \arrow[r] & 1 \text{ }.
\end{tikzcd}
\]
Then, we add the following two arrows
\[
\begin{tikzcd}
a_1:n+2 \arrow[r] &  j  & a_2 : i \arrow[r] & n+1 \text{ }.
\end{tikzcd}
\]
Finally, we get the following quiver,
\[
\begin{tikzcd}
\hat{Q}_{\beta_1 \beta_2} : \quad n+2 \arrow[r] & j \arrow[r] &  \cdots \arrow[r]  & 1 \text{ }, & i \arrow[r] & n+1 \arrow[r] & \cdots \arrow[r] & \arrow[lll, bend right] \text{ }.
\end{tikzcd}
\]

They represent two disjoint cycles in $S_{n+1}$ by replacing $n+2$ by $1$. Denote by $\beta_{1}$ and $\beta_{2}$, $1 \in \beta_1$. So, $\beta=\beta_{1}\beta_{2}\alpha_{2}\dots\alpha_{k}$ is the permutation in $S_{n+1}$ constructed by cutting the arrow $a$.
 	
\item{\textbf{Case 2} (\emph{Join Case}), $a \notin \hat{Q}_{\alpha_{1}}$}
 	
Without loss of generality, we assume that $a \in \hat{Q}_{\alpha_{2}}$. The corresponding quiver for $\alpha_{1}$ and $\alpha_{2}$ are
\[
\begin{tikzcd}
\hat{Q}_{\alpha_1 \alpha_2}: \quad n+1 \arrow[r] &  \cdots \arrow[r]  & 1 \text{ }, & i \arrow[r] & j \arrow[r] & \cdots \arrow[r] & \arrow[lll,bend right] \text{ }.
\end{tikzcd}
\]

Similar to \textbf{Case 1}, we cut the arrow $i \to j$ and we get
\[
\begin{tikzcd}
n+1 \arrow[r] &  \cdots \arrow[r]  & 1 \text{ }, &  j \arrow[r] & \cdots \arrow[r] & i \text{ }.
\end{tikzcd}
\]

Then, we add the following two arrows
\[
\begin{tikzcd}
a_1:n+2 \arrow[r] &  j  & a_2 : i \arrow[r] & n+1 \text{ }.
\end{tikzcd}
\]
Finally, we get the chain
\[
\begin{tikzcd}
\hat{Q}_{\beta_1}: \quad n+2 \arrow[r] &  j \arrow[r] & \cdots \arrow[r]  &  i \arrow[r] & n+1 \arrow[r] & \cdots \arrow[r] & 1 \text{ }.
\end{tikzcd}
\]
It represents a cycle in $S_{n+1}$ by replacing $n+2$ by $1$. Denote by $\beta_{1}$. So, $\beta=\beta_{1}\alpha_{3}\dots\alpha_{k}$ is the permutation in $S_{n+1}$ constructed in this case.
\end{itemize}
\end{con}

The following theorem is a direct result from Construction \ref{202}.
\begin{thm}\label{203}
	For any $\alpha \in S_{n}$, Construction \ref{202} gives $n+1$ distinct permutations in $S_{n+1}$. In fact, if we do it for all $\alpha \in S_{n}$, we will get $(n+1)!$ elements, which are exactly all permutations in the group $S_{n+1}$.
\end{thm}

We define the following notation $[\alpha,j]$, which will be used in the next section.
\begin{defn}\label{204}
	Let $\alpha$ be a permutation in $S_n$. Denote by $[\alpha,j]$ the permutation in $S_{n+1}$ constructed from $\alpha$, where $j$ is an integer, $0 \leq j \leq n$.
	\begin{enumerate}
		\item{$j=0$} $[\alpha,0]$ corresponds to the \textbf{Case 0} in Construction \ref{202}.
		\item{$j>0$} $[\alpha,j]$ corresponds to \textbf{Case 1,2} by cutting the arrow $a$ such that $t(a)=j$.
	\end{enumerate}
\end{defn}

\section{$W$-Operator}
We first review the definition and some properties about the $W$-operator $W([n])$. Details can be found in \cite{MR2864467, Mir}. Then, we will show how to calculate $W([n+1])$ from $W([n])$ by induction, which relates to Construction \ref{202} in the last section. In other words, the structure of the $W$-operator $W([n])$ is closely related to the permutation group $S_n$ and the quivers we studied in \S 2.

Let $X$ be an infinite matrix with variable $X_{ab}$ in the $(a,b)$-entry, i.e. $X:=(X_{ab})_{a\geq1,b\geq1}$. Given a positive integer $k$, denote by $p_{k}$ the trace of $X^{k}$,
\begin{align*}
p_k=\sum_{a_1,\dots,a_k \geq 1} X_{a_1 a_k}X_{a_k a_{k-1}}\dots X_{a_2 a_1}.
\end{align*}
Clearly, $p_{k}$ is a power series in $\mathbb{C}[[X_{ab}]]_{a,b \geq 1}$.

The operator matrix $D$ is the infinite matrix with $D_{ab}$ in the $(a,b)$-entry, where
\begin{align*}
D_{ab}=\sum\limits_{c=1}^{\infty} X_{ac} \frac{\partial}{\partial X_{bc}}.
\end{align*}
In the rest of the paper, we prefer to write $D_{ab}=X_{ac} \frac{\partial}{\partial X_{bc}}$ with the sum over $c$ implied. As differential operators, the normal ordered product of $D_{ab}$ and $D_{cd}$ is defined as
\begin{align*}
:D_{ab}D_{cd}:=X_{ae_{1}}X_{ce_{2}} \frac{\partial}{\partial X_{be_{1}}} \frac{\partial}{\partial X_{de_{2}}}.
\end{align*}
It means that we always calculate the differentiation first. The normal ordered product $:D_{a_{n+1} a_n}D_{a_n a_{n-1}}\dots D_{a_2 a_1}:$ is defined similarly.
\begin{defn}\label{301}
	For any positive integer $n$, the $W$-operator $W([n])$ is defined as follows
	\begin{align*}
	W([n]):=\frac{1}{n}:tr(D^{n}):=\frac{1}{n}\sum_{a_1,\dots,a_n \geq 1} :D_{a_1 a_n}D_{a_n a_{n-1}}\dots D_{a_2 a_1}:.
	\end{align*}
\end{defn}

Before we calculate $W([n+1])$ by induction, we first review some important formulas.
\begin{lem}[\cite{Mir}]\label{302}
	Let $F(p)$ be an element in  $\mathbb{C}[p_1,p_2,\dots]$. We have
	\begin{equation*}
	D_{ab}F(p)=\sum\limits_{k=1}^{\infty} k(X^{k})_{ab} \frac{\partial F(p)}{\partial p_{k}}.
	\end{equation*}
\end{lem}
As a differential operator, the above formula of $D_{ab}$ can be written as
\begin{align*}
	D_{ab}=\sum\limits_{k=1}^{\infty} k(X^{k})_{ab} \frac{\partial }{\partial p_{k}}.
\end{align*}
\begin{lem}[\cite{Mir}]\label{303}
	\begin{equation*}
		D_{cd}(X^{k})_{ab}=\sum_{j=0}^{k-1} (X^j)_{ad}(X^{k-j})_{cb}.
	\end{equation*}
\end{lem}
In particular, Lemma \ref{303} implies the following formula
\begin{align*}
\sum_{k_j=1}^{\infty}D_{a_{n+1}a_{n}}(X^{k_j})_{a_i a_j}
=&\sum\limits_{k_j=1}^{\infty}\sum_{k_n=0}^{k_j-1}(X^{k_n})_{a_{i}a_{n}}(X^{k_j-k_n})_{a_{n+1}a_j}\\
=&\sum\limits_{k_j=1}^{\infty}\sum_{k_n=1}^{\infty}(X^{k_n})_{a_{i}a_{n}}(X^{k_j})_{a_{n+1}a_j}		
\end{align*}
by setting $a=a_i, b=a_j, c=a_{n+1}, d=a_n$.	

\begin{lem}[\cite{MR2864467,Mir}]\label{304}
	\begin{align*}
	D_{a_{n+2}a_{n+1}}D_{a_{i} a_{j}}&=\sum\limits_{k\geq 1, j \geq 0} ((k+j)(X^{j})_{a_{i}a_{n+1}}(X^{k})_{a_{n+2}a_{j}} \frac{\partial}{\partial p_{k+j}})\\
	& +\sum\limits_{k,j\geq1}(kj(X^{k})_{a_{n+2}a_{n+1}}(X^{j})_{a_{i}a_{j}} \frac{\partial ^{2}}{\partial p_{k} \partial p_{j}}),\\
	: D_{a_{n+2}a_{n+1}}D_{a_{i} a_{j}} :&=\sum\limits_{k,j\geq 1} ((k+j)(X^{j})_{a_{i}a_{n+1}}(X^{k})_{a_{n+2}a_{j}} \frac{\partial}{\partial p_{k+j}})\\
	& +\sum\limits_{k,j\geq1}(kj(X^{k})_{a_{n+2}a_{n+1}}(X^{j})_{a_{i}a_{j}} \frac{\partial ^{2}}{\partial p_{k} \partial p_{j}}).
	\end{align*}
\end{lem}

\begin{rem}\label{305}
	The calculation of $D_{a_{n+2}a_{n+1}}D_{a_{i} a_{j}}$ in Lemma \ref{304} is a result of Lemma \ref{302} and \ref{303}. The operator $D_{a_{n+2}a_{n+1}}$ acts on the formula in Lemma \ref{302}. Since $D_{a_{n+2}a_{n+1}}D_{a_{i} a_{j}}$ is a differential operator, the chain rule gives us the result, where the first line comes from the action on the polynomial part (Lemma \ref{303}) and the action on the differential part gives us the second line (Lemma \ref{302}).
	
	Notice that the only difference between the above two formulas in Lemma \ref{304} is that the subscript $j$ in the first sum starts from $0$ in the first formula, while it goes from $1$ in the second one. In fact, the formula of normal ordered product $:D_{a_{n+2}a_{n+1}}D_{a_{i} a_{j}}:$ in Lemma \ref{304} comes from  $D_{a_{n+2}a_{n+1}}D_{a_{i} a_{j}}$ by subtracting an extra term. More precisely, by definition of normal ordered product, we have the following equation
	\begin{align*}
	D_{a_{n+2}a_{n+1}}D_{a_{n+1} a_{n-1}}=
	:D_{a_{n+2}a_{n+1}}D_{a_{n+1} a_{n-1}}:+X_{a_{n+2}e_1}[\frac{\partial}{\partial X_{a_{n+1} e_1}},X_{a_{n+1} e_2}]\frac{\partial}{\partial X_{a_{n} e_2}}.
	\end{align*}
The last term is exactly the term missing in the normal ordered product.

As for the normal ordered product $:D_{a_{n+2} a_{n+1}}\dots D_{a_2 a_1}:$, the reader can use the same approach to calculate $:D_{a_{n+2} a_{n+1}}\dots D_{a_2 a_1}:$ from the product $D_{a_{n+2} a_{n+1}}\dots D_{a_2 a_1}$. In the formula of normal ordered product $:D_{a_{n+2} a_{n+1}}\dots D_{a_2 a_1}:$, the sum always goes from $1$ to infinity, while some subscripts start from $0$ in the formula of $D_{a_{n+2} a_{n+1}}\dots D_{a_2 a_1}$.
\end{rem}

\begin{rem}\label{306}
	In this remark, we will briefly explain how the above formulas relates to Construction \ref{202}. We fix a permutation $\alpha \in S_{n}$, positive integers $k$, $k_j$, and $a_j$ $1 \leq j \leq n$. Let $\hat{Q}_{\alpha}=(\hat{V}_{\alpha},\hat{A}_{\alpha})$ be the quiver we defined in Construction \ref{201}. We consider a special differential operator
\begin{align*}
\prod_{b \in \hat{A}_\alpha}(X^{k_j})_{a_{s(b)}a_{t(b)}} \frac{\partial}{\partial p_k},
\end{align*}
where the polynomial part $\prod_{b \in \hat{A}_\alpha}(X^{k_j})_{a_{s(b)}a_{t(b)}}$ corresponds to the quiver $\hat{Q}_\alpha$.

We want to calculate $D_{a_{n+2} a_{n+1}}\left( \prod_{b \in \hat{A}_\alpha}^{n}(X^{k_j})_{a_{s(b)}a_{t(b)}} \frac{\partial}{\partial p_k}\right)$. By the chain rule, the operator $D_{a_{n+2} a_{n+1}}$ will act on the polynomial part $\prod_{b \in \hat{A}_\alpha}^{n}(X^{k_j})_{a_{s(b)}a_{t(b)}}$ and the differential part $\frac{\partial}{\partial p_k}$ separately.
	\begin{enumerate}
		\item When $D_{a_{n+2} a_{n+1}}$ acts on the differential part, the formula in Lemma \ref{302} gives us the result. In the language of quivers, we add one more arrow $a_{n+2} \rightarrow a_{n+1}$ to the quiver $\hat{Q}_\alpha$, which corresponds to \textbf{Case 0} in Construction \ref{202}.
		\item When $D_{a_{n+2} a_{n+1}}$ acts on the polynomial part, assume that it acts on $(X^{k_j})_{a_i a_j}$. Lemma \ref{303} gives us the calculation. In the language of quivers, it means that we cut the arrow $i \rightarrow j$ and add two arrows $i \rightarrow n$ and $n+1 \rightarrow j$, which corresponds to \textbf{Case 1, 2} in Construction \ref{202}.
	\end{enumerate}
This is the reason why we do Construction \ref{202} in \S 2.
\end{rem}

Now we are ready to calculate the $W$-operator $W([n])$ by induction.

\begin{thm}\label{307}
	$W([n])$ is a well-defined operator on $\mathbb{C}[p_{1},p_{2},\dots]$. It can be written as the sum of $n!$ summations $FS_\alpha$, each of which corresponds to a unique quiver $\hat{Q}_\alpha$, or equivalently a unique permutation $\alpha \in S_n$.
\end{thm}

\begin{proof}
To calculate $W([n])$, we have to figure out the formula of $:D_{a_1 a_{n}}D_{a_{n}a_{n-1}}\dots D_{a_2 a_1}:$ for any positive integer $a_i \geq 1$, $1 \leq i \leq n$. By Lemma \ref{304} and Remark \ref{305}, it is equivalent for us to calculate the product $D_{a_1 a_{n}}D_{a_{n}a_{n-1}}\dots D_{a_2 a_1}$. Since we want to use induction to calculate this product, we replace $D_{a_1 a_{n}}$ by $D_{a_{n+1} a_{n}}$.

Let's start from the base step. When $n=1$, by Lemma \ref{302}, we have
\begin{align*}
D_{a_2 a_1}=\sum\limits_{k_1=1}^{\infty} k_1(X^{k_1})_{a_2 a_1} \frac{\partial}{\partial p_{k_1}}.
\end{align*}
We associate this summation to the quiver
\[
\begin{tikzcd}
\hat{Q}_{(1)}: \quad 2 \arrow[r] &  1 \text{ },
\end{tikzcd}
\]
which corresponds to the subscript of $(X^{k_1})_{a_2 a_1}$. Note that there is only one summation. Thus we define
\begin{align*}
FS'_{(1)}=\sum\limits_{k_1=1}^{\infty} k_1(X^{k_1})_{a_2 a_1} \frac{\partial}{\partial p_{k_1}}.
\end{align*}
Replacing $a_2$ by $a_1$ and taking the sum over $a_1$, we have
\begin{align*}
W([1])=\underbrace{\sum_{k_1 \geq 1} k_1 p_{k_1} \frac{\partial}{\partial p_{k_1}}}_{FS_{(1)}}.
\end{align*}
Denote by $FS_{(1)}$ the summation in $W([1])$ corresponding to $FS'_{(1)}$ in $D_{a_2 a_1}$. Note that $FS'_{(1)}=D_{a_2 a_1}$ in the base case. This is for the base case $n=1$.

We do the next step $n=2$ in fully details. The induction step follows directly from the discussion for $n=2$. When $n=2$, we have to calculate $D_{a_3 a_2} D_{a_2 a_1}$,
\begin{align*}
D_{a_3 a_2} D_{a_2 a_1}=\sum\limits_{k_1=1}^{\infty} \left(D_{a_3 a_2}(k_1(X^{k_1})_{a_2 a_1})\right) \frac{\partial}{\partial p_{k_1}}
+\sum\limits_{k_1=1}^{\infty} k_1(X^{k_1})_{a_2 a_1} \left(D_{a_3 a_2} \circ \frac{\partial}{\partial p_{k_1}}  \right).
\end{align*}
By Lemma \ref{304}, we have
\begin{align*}
D_{a_3 a_2}D_{a_2 a_1}&=\sum\limits_{k_1\geq1,k_2\geq0} ((k_1+k_2)(X^{k_2})_{a_2 a_2}(X^{k_1})_{a_3 a_1} \frac{\partial}{\partial p_{k_1+k_2}})\\
& +\sum\limits_{k_1,k_2\geq1}(k_1 k_2(X^{k_2})_{a_3 a_2}(X^{k_1})_{a_2 a_1}) \frac{\partial ^{2}}{\partial p_{k_1} \partial p_{k_2}}).
\end{align*}
We associate the first summation to the quiver $\hat{Q}_{(1)(2)}$
\[
\begin{tikzcd}
\hat{Q}_{(1)(2)} : \quad 2 \arrow[out=30,in=10,loop]  & \text{ }, & 3 \arrow[r] &  1 \text{ },
\end{tikzcd}
\]
which comes from the subscripts of the polynomial part $(X^{k_2})_{a_2 a_2}(X^{k_1})_{a_{3} a_1}$. Similarly, the second summation corresponds to the quiver $\hat{Q}_{(12)}$
\[
\begin{tikzcd}
\hat{Q}_{(12)} : \quad 3 \arrow[r]  & 2 \arrow[r] &  1 \text{ }.
\end{tikzcd}
\]

We know that $D_{a_3 a_2}$ acting on $(X^{k_1})_{a_2 a_1}$ gives the first summation, which corresponds to \textbf{Case 1} of cutting the arrow $2 \rightarrow 1$ in $\hat{Q}_{(1)}$ in Construction \ref{202}. The same argument holds for the second summation, where $D_{a_3 a_2}$ acts on $\frac{\partial}{\partial p_{k_1}}$ and it corresponds to the \textbf{Case 0} in Construction \ref{202}.

By Lemma \ref{304} and Remark \ref{305}, we know that $:D_{a_3 a_2} D_{a_2 a_1}:$ and $D_{a_3 a_2} D_{a_2 a_1}$ are almost the same and the only difference comes from the term with subscript $j=0$ in the first summation. Hence, we can use quivers to describe the summations of $:D_{a_3 a_2} D_{a_2 a_1}:$ in the same way as $D_{a_3 a_2} D_{a_2 a_1}$. We use the notation $FS'_\alpha$ for the summation corresponding to $\alpha \in S_2$. We have
\begin{align*}
:D_{a_3 a_2}D_{a_2 a_1}:=\sum_{\alpha \in S_2}FS'_\alpha.
\end{align*}
In conclusion, we find that $:D_{a_3 a_2}D_{a_2 a_1}:$ can be written as the sum of two summations, which correspond to quivers $\hat{Q}_\alpha$, $\alpha \in S_2$,
\begin{align*}
:D_{a_3 a_2}D_{a_2 a_1}:&=\underbrace{\sum\limits_{k_1,k_2 \geq1} ((k_1+k_2)(X^{k_2})_{a_2 a_2}(X^{k_1})_{a_3 a_1} \frac{\partial}{\partial p_{k+j}})}_{FS'_{(1)(2)}}\\
& +\underbrace{\sum\limits_{k_1,k_2\geq1}(k_1 k_2(X^{k_2})_{a_3 a_2}(X^{k_1})_{a_2 a_1}) \frac{\partial ^{2}}{\partial p_{k_1} \partial p_{k_2}})}_{FS'_{(12)}}.
\end{align*}
Replacing $a_3$ by $a_1$ and taking the sum over $a_1,a_2$, we have
\begin{align*}
W([2])=\underbrace{\frac{1}{2}\sum_{k_1,k_2 \geq 1}(k_1+k_2)p_{k_1} p_{k_2}\frac{\partial}{\partial p_{k_1+k_2}}}_{FS_{(1)(2)}}+ \underbrace{\frac{1}{2}\sum_{k_1,k_2 \geq 1}k_1 k_2 p_{k_1+k_2}\frac{\partial^2}{\partial p_{k_1} \partial p_{k_2}}}_{FS_{(12)}}.
\end{align*}
We define $FS_{\alpha}$ to be summation of $W([2])$ which corresponds to $FS'_{\alpha}$ in the formula of $:D_{a_3 a_2}D_{a_2 a_1}:$.

When $n=3$, we can calculate the normal ordered product $:D_{a_4 a_3}D_{a_3 a_2}D_{a_2 a_1}:$ from the product $D_{a_4 a_3}D_{a_3 a_2}D_{a_2 a_1}$ similarly. We use the operator $D_{a_4 a_3}$ acting on $D_{a_3 a_2}D_{a_2 a_1}$. Since the polynomial part of each summation is a product of two terms (or the corresponding quiver has two arrows), we will get three new summations by the chain rule where two comes from the polynomial part and one from the differential part. By Theorem \ref{203} and Remark \ref{306}, each of the new summation corresponds to a unique permutation in $S_3$ and this correspondence is a one-to-one correspondence. Thus the normal ordered product $:D_{a_4 a_3}D_{a_3 a_2}D_{a_2 a_1}:$ can be written as the sum of summations $FS'_{\alpha}$ labelled by permutations $\alpha$ in $S_3$. Replacing $a_4$ by $a_1$ and taking the sum over $a_1,a_2,a_3$, we get the formula for $W([3])$. We define $FS_\alpha$ to be the summation of $W([3])$, which corresponds to the summation $FS'_\alpha$. In this way, the operator $W([n])$ can be written as the sum of $n!$ summations by induction, and each summation corresponds to a unique permutation in $S_n$.

At the end of the proof, we will show that $W([n])$ is a well defined operator on $\mathbb{C}[p_1,p_2,\dots]$ by giving two examples. When $n=1$, consider $D_{a_2 a_1}$. Let $a_2=a_1$ and take the sum over $a_1$. We have
\begin{align*}
	W([1])=\sum_{k_1} k_1 p_{k_1} \frac{\partial}{\partial p_{k_1}}.
\end{align*}
Now we consider $:D_{a_3 a_2}D_{a_2 a_1}:$. Let $a_3=a_1$ and sum over $a_1,a_2$. We have
\begin{equation*}
W([2]) = \frac{1}{2}\sum_{a_1,a_2 \geq 1}:D_{a_1 a_2}D_{a_2 a_1}:
=\frac{1}{2}\sum_{k_1,k_2\geq1}((k_1+k_2)p_{k_1} p_{k_2} \frac{\partial}{\partial p_{k_1+k_2}}+k_1 k_2 p_{k_1+k_2}\frac{\partial^2}{\partial p_{k_1} \partial p_{k_2}}).
\end{equation*}
Clearly, $W([1])$ and $W([2])$ are well-defined operator on $\mathbb{C}[p_1,p_2,\dots]$. The operator $W([n])$ can be proved to be a well-defined operator on $\mathbb{C}[p_1,p_2,\dots]$ by induction.
\end{proof}

\section{Degree of Summation $FS_\alpha$}

\begin{defn}\label{401}
	Given any summation $FS_\alpha$ of $W([n])$, define $dP(FS_\alpha)$ to be the degree of its polynomial part and $dD(FS_\alpha)$ the order of its differential part. The degree $d(FS_\alpha)$  of the summation $FS_\alpha$ is $d(FS_\alpha)=dP(FS_\alpha)+dD(FS_\alpha)$.
\end{defn}
We give an easy example to explain this definition. Consider the summation \begin{align*}
FS_{(1)(2)}=\frac{1}{2}\sum_{k_1,k_2 \geq 1}(k_1+k_2)p_{k_1} p_{k_2}\frac{\partial}{\partial p_{k_1+k_2}}.
\end{align*}
We have
\begin{align*}
dP(FS_\alpha)=2, \quad dD(FS_\alpha)=1, \quad d(FS_\alpha)=3.
\end{align*}
Similarly, the degree data of $FS_{(12)}$ is
\begin{align*}
dP(FS_\alpha)=1, \quad dD(FS_\alpha)=2, \quad d(FS_\alpha)=3.
\end{align*}

The following lemma describes the relation between the degree of $FS_\beta$ and $FS_\alpha$, where $\beta=[\alpha,j]$ (see Definition \ref{204}).
\begin{lem}\label{402}
	For any $\alpha \in S_{n}$,
	\begin{enumerate}
		\item If $[\beta] = [\alpha,0]$, we have
		\begin{align*}
		dP(FS_\beta)=dP(FS_\alpha), \quad dD(FS_\beta)=dD(FS_\alpha)+1, \quad d(FS_\beta)=d(FS_\alpha)+1.
		\end{align*}
		\item If $[\beta] = [\alpha,j]$ and $j$ is a vertex in the chain of $\hat{Q}_\alpha$, then, we have
		\begin{align*}
		dP(FS_\beta)=dP(FS_\alpha)+1, \quad dD(FS_\beta)=dD(FS_\alpha),\quad d(FS_\beta)=d(FS_\alpha)+1.
		\end{align*}
		\item If $[\beta] = [\alpha,j]$ and $j$ is not a vertex in the chain of $\hat{Q}_\alpha$, we have
		\begin{align*}
		dP(FS_\beta)=dP(FS_\alpha)-1,\quad dD(FS_\beta)=dD(FS_\alpha),\quad d(FS_\beta)=d(FS_\alpha)-1.
		\end{align*}
	\end{enumerate}
\end{lem}

\begin{proof}
	Notice that the $dP(FS_\alpha)$ is exactly the number of disjoint cycles of $\alpha$.
	
	When $j=0$, the differential degree of $FS'_\beta$ increases by one by Lemma \ref{302}. The disjoint cycle of $\beta=[\alpha,0]$ is the same as $\alpha$. Hence, $dP(FS_\beta)=dP(FS_\alpha)$.
	
	When $j \geq 1$, Lemma \ref{303} and Remark \ref{306} implies that the operator $D_{a_{n+n}a_{n+1}}$ fixes the differential degree. Now we consider the polynomial degree. If $j$ is in the chain of $\hat{Q}_\alpha$, \textbf{Case 1} in Construction \ref{202} tells us that $\beta$ has one more disjoint cycle than $\alpha$. When $j$ is not in the chain of $\hat{Q}_\alpha$, $\beta$ corresponds to \textbf{Case 2} and $dP(FS_\beta)=dP(FS_\alpha)-1$. This finishes the proof of the lemma.
\end{proof}

From the above lemma, it is easy to find that the highest degree of summations in $W([n])$ is $n+1$ and the other possible degrees are $n-1,n-3,\dots \text{ }.$

\begin{defn}[Ordinary Summation]\label{403}
	Let $\alpha$ be a permutation in $S_n$. We say that $FS_\alpha$ is an \emph{ordinary summation} (OS), if $d(FS_\alpha)=n+1$.
	
	An ordinary summation $FS_\alpha$ is of \emph{type $(r,n-r+1)$} if $dP(FS_\alpha)=r$ and $dD(FS_\alpha)=n-r+1$.
\end{defn}

\section{Noncrossing Permutation}
In this section, we prove that if the permutation $\alpha$ is a noncrossing permutation if and only if $FS_\alpha$ is of maximal degree. As a corollary of this correspondence, we show that the number of ordinary summations with maximal degree is the Catalan number.

The noncrossing permutation comes from noncrossing partition with respect to a fixed order of objects. Recall that a partition of $[n]=\{1,\dots,n\}$ is \emph{noncrossing} if whenever four elements, $1 \leq a < b<c<d \leq n$, are such that $a,c$ are in the same block and $b,d$ are in the same block, then the two blocks coincide. With respect to the natural order of integers, each noncrossing partition corresponds to a unique permutation, where each block corresponds to a disjoint cycle and the order $i > j$ implies an arrow $i \rightarrow j$ in the disjoint cycle. Here is the definition of noncrossing permutation.

\begin{defn}[Noncrossing permutation]\label{501}
	Let $\alpha$ be a permutation in $S_n$. Let $\alpha=\alpha_1 \dots\alpha_r$ be the decomposition of $\alpha$ into disjoint cycles. The permutation $\alpha$ is a \emph{noncrossing permutation} if it  satisfies the following conditions
	\begin{enumerate}
		\item[$(*_1)$] For each arrow $a$ in the unique chain of $\hat{Q}_\alpha$, we have $t(a)<s(a)$, and there is only one arrow $b$ in each loop of $\hat{Q}_\alpha$ such that $s(b)<t(b)$.
		\item[$(*_2)$] Given any two distinct cycles $\alpha_{i}$ and $\alpha_{j}$, they satisfy at least one of the following conditions:
		\begin{enumerate}
			\item pick an arbitrary element $m$ in $\alpha_{i}$, either $m>n$ for any $n$ in $\alpha_{j}$ or $m<n$ for any $n$ in $\alpha_j$;
			\item pick an arbitrary element $m$ in $\alpha_{j}$, either $m>n$ for any $n$ in $\alpha_{i}$ or $m<n$ for any $n$ in $\alpha_i$.
		\end{enumerate}	
	\end{enumerate}
\end{defn}

Condition $(*_1)$ means that we have an order on the finite set which determines the permutation. Condition $(*_2)$ corresponds to the \emph{noncrossing} condition. If $\alpha_i$ and $\alpha_j$ satisfy only one of the condition $(a),(b)$, then then one is contained in the other one. If they satisfy both these two conditions, they are disjoint. For instance, consider the following examples,
\begin{align*}
\tau_1=(123)(45), \quad \tau_2=(125)(34), \quad \tau_3=(124)(35).
\end{align*}
The two disjoint cycles in $\tau_1$ satisfies both these two conditions. If we consider the integers are in a straight line
\begin{align*}
( \quad 5 \quad 4 \quad ) \quad (\quad 3 \quad 2 \quad 1 \quad ),
\end{align*}
the disjoint cycles $(123)$ and $(45)$ are disjoint.

The permutation $\tau_2$ is an example that $(125)$ and $(34)$ only satisfy one of the conditions $(a),(b)$. It is obvious that the disjoint cycle $(34)$ is contained in $(125)$,
\begin{align*}
( \quad 5 \quad ( \quad 4 \quad 3 \quad ) \quad 2 \quad 1 \quad ).
\end{align*}

Now we will prove that $FS_\alpha$ is an ordinary summation if and only if $\alpha$ is a noncrossing permutation in several steps.

\begin{lem}\label{502}
	Given $\alpha \in S_n$, if $FS_\alpha$ is an OS, then $\alpha$ satisfies the condition $(*_1)$.
\end{lem}

\begin{proof}
	We prove this lemma by induction on $n$. For the base step $n=1$, $\hat{Q}_{(1)}$ is the only quiver and $FS_{(1)}$ is an OS. There is only one arrow $2 \rightarrow 1$ in the quiver $\hat{Q}_{(1)}$. Clearly, $(1)$ satisfies the condition $(*_1)$.
	
	By induction, we assume that, for all $\alpha \in S_{k-1}$, if $FS_\alpha$ is an OS, then $\alpha$ satisfies $(*_1)$. Let $\beta=[\alpha,j]$ be a permutation in $S_k$. The condition $FS_\beta$ is an OS implies that $FS_\alpha$ is also an OS. Indeed, if $FS_\alpha$ is not an OS, then $d(FS_\alpha)<k$. By Lemma \ref{402}, $d(FS_\beta)<k+1$, violating the fact that $\beta$ is an OS.
	
	Let $\alpha = \alpha_{1} \dots \alpha_{r}$ be the decomposition of $\alpha$ into disjoint cycles with $1 \in \alpha_{1}$. By Lemma \ref{402}, the integer $j$ is either zero or the target of some arrow in the chain of $\hat{Q}_\alpha$. Now we discuss these two cases.
	\begin{enumerate}		
		\item	If $j=0$, then $\beta=\beta_{1}\alpha_{2} \dots\alpha_{r}$, where $\hat{Q}_{\beta_{1}}$ is constructed from $\hat{Q}_{\alpha_{1}}$ by adding another arrow $k+1 \to k$. By induction, the statement is true.
		
		\item	If $j \neq 1$, then $\beta$ is constructed from $\alpha$ by cutting the arrow $a:i \to j$, which is an arrow $a$ in the chain of $\hat{Q}_\alpha$. We use the same notation as \textbf{Case 1} in Construction \ref{202}. Let $\beta=\beta_{1}\beta_{2}\alpha_{2}\dots \alpha_{r}$. The quiver $\hat{Q}_{\beta_{1}}$ of the cycle $\beta_1$ is
		\[
		\begin{tikzcd}
		\hat{Q}_{\beta_1} : \quad k+2 \arrow[r]  & j \arrow[r] & \cdots \arrow[r] &  1 \text{ },
		\end{tikzcd}
		\]
		where $j \to \dots \to 1$ is a subquiver of $\alpha_{1}$. Hence, all arrows in this chain satisfy that the source is larger than the target. The quiver $\hat{Q}_{\beta_2}$ is
		\[
		\begin{tikzcd}
		\hat{Q}_{\beta_2} : \quad i \arrow[r]  & k+1 \arrow[r] & \cdots \arrow[r] &  \arrow[lll, bend right] \text{ },
		\end{tikzcd}
		\]
		where $k+1 \rightarrow\dots \rightarrow i$ is a subquiver of $\alpha_1$ by construction. So the only arrow $a$ in the cycle $\hat{Q}_{\beta_2}$ satisfying $s(a)<t(a)$ is $i \rightarrow k+1$. Hence, the statement is true for $n=k$.
	\end{enumerate}
\end{proof}

\begin{lem}\label{503}
	If $FS_\alpha$ is an OS, then $\alpha$ satisfies the condition $(*_2)$.
\end{lem}

\begin{proof}
	Similar to the proof of Lemma \ref{502}, we prove this lemma by induction on the permutation group $S_n$. When $n=1$, it is clear that the unique permutation $(1)$ in $S_{1}$ satisfies the condition $(*_2)$.
	
	By induction, we assume that for all $\alpha \in S_{k-1}$ if $FS_\alpha$ is an OS, then $\alpha$ satisfies $(*_2)$. Let $\beta=[\alpha,j] \in S_k$. Let $\alpha = \alpha_{1} \dots\alpha_{r}$ be the decomposition of $\alpha$ into disjoint cycles. We will prove that if $FS_\beta$ is an OS, then $\beta$ satisfies the condition $(*_2)$. Before we give the proof, we remind the reader that if $[\beta]=[\alpha,j]$ and $FS_{\beta}$ is an OS, then $FS_{\alpha}$ is also an OS. This property comes from the proof of Lemma \ref{502}.
	
	If $j=0$, then $\beta=\beta_{1}\alpha_{2}\dots\alpha_{r}$, where $\hat{Q}_{\beta_{1}}$ is constructed from $\hat{Q}_{\alpha_{1}}$ by adding another arrow $k+1 \to k$. In other words, we put another element $k$ into the cycle $\alpha_{1}$ (see Construction \ref{202}). By assumption that any two disjoint cycles of $\alpha \in S_{k-1}$ satisfy at least one of the conditions, we only have to check whether the pair $(\beta_{1},\alpha_{i})$ satisfies the condition $(*_2)$, $2 \leq i \leq r$. Since $\alpha_1$ contains the smallest element $1$, so if $\alpha_1$ and $\alpha_i$ are disjoint, then any element in $\alpha_1$ is smaller than any element in $\alpha_i$. Since $k$ is the largest element, so the statement is true for $\beta_1$ and $\alpha_i$. Now we consider that $\alpha_1$ and $\alpha_i$ are not disjoint. Since $1$ is contained in $\alpha_1$, so $\alpha_i$ is contained in $\alpha_1$. Clearly, it still holds for $\beta_1$ and $\alpha_i$. So, $(\beta_{1},\alpha_{i})$ satisfies the condition $(*_2)$.
	
	Now let's consider the case that $\beta$ is constructed from $\alpha$ by cutting the arrow $a:i \to j$ lying in the chain of $\alpha$. We use the same notation as \textbf{Case 1} in Construction \ref{202}. Let $\beta=\beta_{1}\beta_{2}\alpha_{2}\dots\alpha_{r}$. So, we have to check whether the following three types of pairs satisfy the condition:
	\begin{align*}
	(\beta_{1},\beta_{2}), \quad (\beta_{1},\alpha_{i}), \quad (\beta_{2},\alpha_{i}),
	\end{align*}
	where $ 2 \leq i \leq r$.
	\begin{itemize}
		\item $(\beta_{1},\beta_{2})$
		
		Since $FS_\alpha$ is OS, all arrows $a$ in $\hat{Q}_{\alpha_{1}}$ satisfy $t(a)<s(a)$. Hence, when cutting the arrow $i \to j$, any elements in $\beta_{2}$ is larger than any elements in $\beta_{1}$. It is true in this case.
		
		\item $(\beta_{1},\alpha_{i})$
		
		By induction, we know that the lemma is true for $(\alpha_{1},\alpha_{i})$, $2 \leq i \leq r$. Since the elements of $\beta_{1}$ is a subset of the elements of $\alpha_{1}$, so it is true for $(\beta_{1},\alpha_{i})$, $2 \leq i \leq r$.
		
		\item $(\beta_{2},\alpha_{i})$
		
		If $\beta_{2}$ is a single disjoint "one cycle" $(k)$, the statement is true. If $\beta_2 \neq (k)$, assume the largest element in $\beta_{2}$ except $k$ is $\phi$. If $\phi$ is smaller than the smallest element in $\alpha_{i}$, then any element $u$ except $k$ in $\beta_2$ $u$ is smaller than any element in $\alpha_i$. Also, $k$ is larger than any element in $\alpha_{i}$. Hence, the statement is true in this case. Now let's consider the case that $\phi$ is larger than the smallest element in $\alpha_{i}$. By construction, $\phi$ is an element in $\alpha_{1}$, which contains $1$. Hence, $\phi$ is larger than any elements in $\alpha_{i}$ by induction. Similarly, any other elements in $\beta_2$ is larger or smaller to all elements in $\alpha_i$ by induction. So, the statement is true.
	\end{itemize}
	This finishes the proof of this lemma.
\end{proof}

\begin{thm}\label{504}
	The summation $FS_{\alpha}$ is an OS if and only if $\alpha$ is a noncrossing permutation.
\end{thm}

\begin{proof}
	The "only if" part is exactly Lemma \ref{502} and \ref{503}. So, we only have to prove the "if" part. We prove this theorem by induction on $n$.
	
	When $n=1$, it is easy to prove, since $(1)$ is the only permutation.
	
	By induction, we assume that if $\alpha \in S_{k-1}$ satisfies the condition $(*)$, then $FS_{\alpha}$ is an OS. We will prove that if $\beta \in S_{k}$ satisfies the condition $(*)$, then $FS_\beta$ is an OS. Assume $[\beta]=[\alpha,j]$ for some $\alpha$ in $S_{n-1}$ and some nonnegative integer $j$. We claim that $j$ is $0$ or in the chain of $\hat{Q}_\alpha$ (\textbf{Claim 1}). Also, we claim that $\alpha$ is a noncrossing permutation (\textbf{Claim 2}).
	
	By \textbf{Claim 1}, $j$ is $0$ or in the chain of $\hat{Q}_\alpha$. By \textbf{Claim 2}, $FS_\alpha$ is an ordinary summation. By Construction \ref{202} and Lemma \ref{402}, we know $FS_\beta$ is also an OS. Now we are going to prove these two claims.
	
	\begin{proof}[Proof of \textbf{Claim 1}]
		If not, $\beta$ is constructed from $\alpha$ by cutting arrow $a:i \to j$ which is not in the chain of $\hat{Q}_\alpha$. Hence, by \textbf{Case 2} in Construction \ref{202}, we will get a long chain
		\[
		\begin{tikzcd}
		k+1 \arrow[r]  & j \arrow[r] & \cdots \arrow[r] &  i \arrow[r] & k \arrow[r] & \cdots \arrow[r] & 1 \text{ }.
		\end{tikzcd}
		\]
		
		In this chain, we have $i<k$, which contradicts with our assumptions that $\beta$ satisfies the condition $(*_1)$. So, $j$ must be in the chain of $\hat{Q}_\alpha$ or $j=0$.	
	\end{proof}
	
	\begin{proof}[Proof of \textbf{Claim 2}]
		By \textbf{Claim 1}, we know that $j=0$ or $j$ is in the chain of $\hat{Q}_\alpha$. If $j=0$, it is easy to prove $\alpha$ is a noncrossing permutation. We leave it for the reader. Now we assume that $j$ is in the chain of $\hat{Q}_\alpha$. With the same notation as in Construction \ref{202}, let $\beta=\beta_{1}\beta_{2}\alpha_{2}\dots\alpha_{r}$ with $1 \in \beta_{1}$.
		
		First, we have to check $\alpha$ satisfies the condition $(*_1)$. By the assumption of $\beta$, there is exactly one arrow $a$ in the quiver of $\alpha_i$ such that $t(a)>s(a)$, where $2 \leq i \leq r$. So, we have to show all arrows $a$ in the chain of $\hat{Q}_\alpha$ satisfying $t(a)<s(a)$. We assume that there is an arrow $a$ in the chain of $\hat{Q}_\alpha$ such that $s(a)<t(a)$. If $t(a) \neq j$, then this arrow will be in either $\beta_{1}$ or $\beta_{2}$, which contradicts with the assumption of $\beta$. If $t(a)=j$, then we get $\beta_{1}$
		\[
		\begin{tikzcd}
		\hat{Q}_{\beta_1}: \quad k+1 \arrow[r]  & j \arrow[r] & \cdots \arrow[r] &   1 \text{ }.
		\end{tikzcd}
		\]
		and $\beta_{2}$
		\[
		\begin{tikzcd}
		\hat{Q}_{\beta_2} : \quad i \arrow[r]  & k \arrow[r] & \cdots \arrow[r] &  \arrow[lll, bend right] \text{ }.
		\end{tikzcd}
		\]
		Since $k>j>i$, so $(\beta_1,\beta_2)$ does not satisfy the condition $(*_2)$. Hence, we have $t(a)<s(a)$ for each arrow $a$ in the chain of $\hat{Q}_\alpha$ and there is exactly one arrow $b$ in each loop of $\hat{Q}_\alpha$ such that $s(b)<t(b)$.
		
		Now, we are going to prove that $\alpha$ satisfies the condition $(*_2)$. The problem pair is $(\alpha_{1},\alpha_{i})$, $2 \leq i \leq r$. By assumption, $\beta_1$ contains the smallest element $1$ and $\beta_2$ contains the element $k$. Hence, by Construction \ref{202} and Lemma \ref{503}, we know that any element in $\beta_1$ is smaller than any element in $\beta_2$. Since $\beta$ is a noncrossing permutation, so for any cycle $\alpha_i$, $2 \leq i \leq r$, there are three possible cases
		\begin{itemize}
			\item $\alpha_i$ is contained in $\beta_1$, i.e. if we pick an arbitrary element $m$ in $\beta_1$, then we have $m>n$ for any $n$ in $\alpha_i$ or $m<n$ for any $n$ in $\alpha_i$;
			\item $\alpha_i$ is contained in $\beta_2$, i.e. if we pick an arbitrary element $m$ in $\beta_2$, then we have $m>n$ for any $n$ in $\alpha_i$ or $m<n$ for any $n$ in $\alpha_i$;
			\item $\alpha_i$ is disjoint with $\beta_1$ and $\beta_2$, i.e. any element in $\alpha_i$ is larger than any element in $\beta_1$ and smaller than any element in $\beta_2$.
		\end{itemize}
		In the first case, if $\alpha_i$ is "contained" in $\beta_1$, then any element in $\beta_2$ is larger than any element in $\alpha_i$, because the element in $\beta_2$ is always larger than the element in $\beta_1$. By the construction of $\alpha_1$, the condition is true for $(\alpha_1,\alpha_i)$. The same argument holds for the second case. For the third case, $\beta_1$ and $\beta_2$ are constructed from $\alpha_1$ by cutting the arrow with target $j$ and add another element $k$. Hence, $\alpha_i$ is "contained" in $\alpha_1$. Hence, $\alpha$ satisfies the condition $(*_2)$.
	\end{proof}	
\end{proof}

\begin{cor}\label{505}
	The number of $(r,n-r+1)$-type OS in $W([n])$ is the Narayana number:
	\begin{align*}
	\frac{1}{n+1}{n+1 \choose r}{n-1 \choose r-1}.
	\end{align*}
	The number of all ordinary summations in $W([n])$ is the Catalan number
	\begin{align*}
	\sum_{r \geq 1}^{n}\frac{1}{n+1}{n+1 \choose r}{n-1 \choose r-1}=\frac{1}{n+1}{2n \choose n}.
	\end{align*}
\end{cor}

\begin{proof}
	By Theorem \ref{504}, there is a one-to-one correspondence between the ordinary summation and noncrossing permutation (also noncrossing partition). The number of $(r,n-r+1)$-type OS is exactly the number of noncrossing partition with $r$ blocks, which is the Narayama number $\frac{1}{n+1}{n+1 \choose r}{n-1 \choose r-1}$ \cite{Sim}. The number of all ordinary summations in $W([n])$ is the Catalan number, which is the sum of Narayama number.
\end{proof}

\textbf{Acknowledgements}\\

This project is based on part of the author's PhD thesis. The author want to thank his advisor Maarten Bergvelt for helpful discussions and invaluable suggestions.

\bigskip
\noindent\small{\textsc{Department of Mathematics, Sun Yat-Sen University}\\
135 Xingang W Rd, BinJiang Lu, Haizhu Qu, Guangzhou Shi, Guangdong Sheng, China}\\
\emph{E-mail address}:  \texttt{sunh66@mail.sysu.edu.cn}
\end{document}